\documentclass[10pt,a4paper,reqno,oneside]{amsart}
\usepackage{config}

\begin{document}
\title[Ping-pong for HNN-extension of free group]{Ping-pong for basis-conjugating HNN-extension of free group}

\author{Vasily Ionin}
\address{St. Petersburg Department \\ of Steklov Institute of Mathematics, \\ Fontanka 27, St. Petersburg 191011, Russia}
\email{ionin.code@gmail.com}

\date{\today}

\begin{abstract}
We isolate a tractable class of HNN-extensions of a free group---namely, multiple HNN-extensions by basis-conjugating embeddings.
For this class, we construct a normal form and establish a practical version of the ping-pong lemma that provides verifiable sufficient conditions for a set of elements to generate a free subgroup.

We then apply these results to the pure braid group \(P_{n+1}\), exploiting its well-known decomposition as a semidirect product of free groups.
Our approach yields new families of free subgroups within the first two factors~\(F_n \rtimes F_{n-1}\) of this decomposition.
\end{abstract}

\thanks{This work was supported by the Ministry of Science and Higher Education of the Russian
Federation (agreement~\texttt{075-15-2025-344} dated 29/04/2025 for Saint Petersburg Leonhard Euler International Mathematical Institute at PDMI RAS).
The work was supported by the Theoretical Physics and Mathematics Advancement Foundation ``BASIS'',
project no.~\texttt{24-7-1-26-3}.}
\maketitle

\section{Introduction}

In this paper, we study the small subclass of (multiple) partial ascending HNN-extensions of finitely generated free groups
that originates from the notion of basis-conjugating automorphisms.
We find sufficient conditions for a set of elements in such a group to generate a free subgroup.
We specialize our results to the pure braid group in order to produce new families of free subgroups.

We begin with a review of basic concepts involved in our work, in order to introduce notation and clarify our motivation.

\subsection{Ascending HNN-extension of free group}
Let~\(F\) be a free group of finite rank.
If~\(\varphi \colon F \to F\) is an automorphism, one may form its \textit{mapping torus}
\[M_\varphi = F \rtimes_\varphi \ZZ = \langle F, t \mid a^t = \varphi(a), \; a \in F\rangle.\]
The group~\(M_\varphi\) has been studied extensively.
Its structure is similar to the fundamental group of a closed three-dimensional manifold fibered over the circle;
see~\cite{FH99}.
If~\(\varphi \colon F \to F\) is a monomorphism,
we have an~\textit{ascending HNN-extension}
\[F *_\varphi = \langle F, t \mid a^t = \varphi(a), \; a \in F \rangle.\]
An ascending HNN-extension of a free group is a common generalization of
the mapping tori of free group automorphisms and the Baumslag--Solitar group
\[BS(1, n) = \langle a, t \mid a^t = a^n \rangle.\]
The ascending HNN-extensions of free groups are of particular importance in the theory of one-relator groups;
see an overview in \cite{Lo25}.
There are some interesting structural results concerning ascending HNN-extensions of free groups.
For example, in \cite{BS05}, it is proved that these groups are residually finite, and
in~\cite{Ka00}, it is shown that most of these groups are word-hyperbolic.

\subsection{Partial ascending HNN-extension of free group}
Let~\(F(Y)\) be a free group of finite rank with a fixed basis~\(Y\).
If~\(Y' \subset Y\) is a subset and~\(\varphi \colon F(Y') \to F(Y)\) is an embedding of groups,
we have a~\textit{partial ascending HNN-extension}
\[\langle F(Y), t \mid y^t = \varphi(y), \; y \in Y' \rangle.\]
This generalization arises naturally in the study of the ordinary ascending HNN-extension of a free group; see \cite{CW24, FH99}.

\subsection{Multiple HNN-extension}
If~\(K\) is a group and~\(\{\varphi_i \mid 1 \le i \le n\}\) is a collection of embeddings of subgroups~\(\{L_i \mid 1 \le i \le n\}\) of \(K\) into \(K\), then the group
\[\langle K, x_1, \ldots, x_n \mid y^{x_i} = \varphi_i(y), \; 1 \le i \le n, \, y \in L_i\rangle\]
is called a \textit{multiple HNN-extension} of \(K\) with \text{stable letters}~\(x_i, \; 1 \le i \le n\).
The structure of subgroups of a multiple HNN-extension is studied in \cite{KS71}.
However, in the context of ascending HNN-extensions of free groups,
the case of multiple extensions has been studied considerably less.
For example, it is unknown whether multiple ascending HNN-extensions of free groups are residually finite; see Problem~4.4 in \cite{BS05}.

\subsection{Basis-conjugating embedding}
Given a subset~\(Y' \subset Y\) and a choice of an element~\(w_y \in F(Y)\) for each~\(y \in Y'\),
we can form a homomorphism
\[\varphi \colon F(Y') \to F(Y), \; \varphi(y) \defeq y^{w_y}, \, y \in Y'.\]
We call such a homomorphism a \textit{basis-conjugating embedding}.
Clearly, this homomorphism is injective, since it induces a monomorphism on the abelianizations.
This notion is related to basis-conjugating automorphisms of free groups.
The subgroup~\(\mathrm{P}\Sigma_{|Y|} \subset \Aut(F(Y))\) consisting of basis-conjugating automorphisms
has been studied extensively in~\cite{CPVW08, JMM06, Mv86}.

\subsection{Overview of results}
Let~\(n \ge 1\) and for every index~\(1 \le i \le n\) we are given a subset
\[Y_i = \{y_{i1}, \ldots, y_{im_i}\} \subset Y\]
and two basis-conjugating embeddings
\[\varphi_i, \psi_i \colon F(Y_i) \to F(Y).\]
Then, for every~\(i\) and~\(1 \le j \le m_i\) there are unique elements~\(w_{ij}, v_{ij} \in F(Y)\)
such that
\[\varphi_i(y_{ij}) = y_{ij}^{w_{ij}}, \; \psi_i(y_{ij}) = y_{ij}^{v_{ij}},\]
and the words~\(y_{ij}^{w_{ij}} = w_{ij}^{-1}y_{ij} w_{ij}\) and~\(y_{ij}^{v_{ij}} = v_{ij}^{-1}y_{ij} v_{ij}\)
are reduced.
We can form a multiple HNN-extension of~\(F(Y)\) by identifying the images of~\(F(Y_i)\) in~\(F(Y)\):
\begin{align}\label{eq:presentation_of_group}
\begin{split}
G &= \langle F(Y), x_1, \ldots, x_n \mid \varphi_i(y)^{x_i} = \psi_i(y), \; 1 \le i \le n, \, y \in Y_i \rangle \\
&= \langle Y, x_1, \ldots, x_n \mid (y_{ij}^{w_{ij}})^{x_i} = y_{ij}^{v_{ij}}, \; 1 \le i \le n, \, 1 \le j \le m_i \rangle.
\end{split}
\end{align}
The basic facts about HNN-extensions say that the natural homomorphism \(F(Y) \to G\) is injective and that the group \(G\) is torsion-free; see \cite{HNN49}.
Note that the abelianization of \(G\) is free abelian of rank \(|Y| + n\).

Let~\(X = \{x_1, \ldots, x_n\}\). We call the set~\(Y^\pm \cup X^\pm\) an \textit{extended alphabet}. Here are our main results.
\begin{Theorem}[Normal form for~\(G\)]
\label{t:diamond_lemma_for_G}
Every element of~\(G\) has a unique reduced expression in the extended alphabet with no subwords of the form
\[x_i v_{ij}^{-1} y_{ij}^{\pm 1}\; \text{and}\; x_i^{-1} w_{ij}^{-1} y_{ij}^{\pm 1}.\]
\end{Theorem}
\begin{Theorem}[Ping-pong lemma for~\(G\)]
\label{t:ping_pong_condition}
Let~\(T_1, \ldots, T_k\) be pairwise disjoint nonempty subsets of~\(X\). If~subgroups~\(A_1, \ldots, A_k \subset G\) satisfy the conditions~\(A_t \subset \langle Y, T_t \rangle\) and~\(A_t \cap \langle Y \rangle = 1\) for every~\(1 \le t \le k\),
then the natural homomorphism
\[A_1 * \ldots * A_k \to G\]
is injective.
\end{Theorem}

Recall that the pure braid group~\(P_n\) is generated by the braids~\(A_{i,j}\) with~\(1 \le i < j \le n\).
There is a~decomposition
\[P_{n+1} = F_n \rtimes F_{n-1} \rtimes \ldots \rtimes \ZZ.\]
The summand~\(F_n \rtimes F_{n-1}\) corresponds to a
subgroup~\(P_{n+1}^{(2)} \subset P_{n+1}\) consisting of braids that become trivial after the last two strands were removed.
We write
\[P_{n+1}^{(2)} = F(x_1, \ldots, x_n) \rtimes F(y_1, \ldots, y_{n-1}),\]
where~\(x_i \defeq A_{i,n+1}\) and~\(y_i \defeq A_{i,n}\).
Note that the abelianization of this group is free abelian of rank \(2n-1\),
generated by the images of~\(x_i\)'s and~\(y_j\)'s. Roughly speaking, the abelianization homomorphism
\[P_{n+1}^{(2)} \to (P_{n+1}^{(2)})_\mathrm{ab}\]
assigns to a pure braid~\(w\),
given by a word in the alphabet~\(Y^\pm \cup X^\pm\), the tuple of exponent sums of the letters.
We denote these sums by~\(\exp_{x_i}(w)\) and~\(\exp_{y_j}(w)\).

The subgroup of~\(P_{n+1}^{(2)}\) generated by the elements~\(\{x_1, \ldots, x_n\}\) is free of rank \(n\).
We use our above results to show that one can replace the generators~\(\{x_1, \ldots, x_{n-1}\}\)
by suitable braids~\(\{w_1, \ldots, w_{n-1}\}\) while keeping the resulting subgroup
\(\langle w_1, \ldots, w_{n-1}, x_n \rangle\)
free.
\begin{Theorem}
\label{t:main_theorem_braids}
Let~\(w_1, w_2, \ldots, w_{n-1}\) be braids from~\(P_{n+1}^{(2)}\)
such that for each~\(1 \le i < n\) one has\textup{:}
\begin{enumerate}
\item~\(\exp_{x_j}(w_i) \ne 0\) if and only if~\(j = i\)\textup{;}
\item~\([w_i, x_n] \ne 1\).
\end{enumerate}
Then the subgroup~\(\langle w_1, \ldots, w_{n-1}, x_n \rangle\) of~\(P_{n+1}\) is free of rank~\(n\).
\end{Theorem}
\begin{Remark}
Condition (2) is crucial. For example, for~\(w_i = y_i x_i\) one has
\[\langle y_1x_1, \ldots, y_{n-1} x_{n-1}, x_n \rangle \cong F_{n-1} \times \ZZ.\]
\end{Remark}
It is easy to show that the elements~\(w_1, \ldots, w_{n-1}\) generate a free group of rank~\(n-1\), and that for every~\(i\) the subgroup~\(\langle w_i, x_n \rangle\) is free of rank~\(2\). We use our method to show that the full subgroup is also free.
We also show that,
despite the fact that the lower central series Lie ring of~\(P_{n+1}^{(2)}\) behaves in a controlled way, our criterion implies the freeness of some subgroups that cannot be deduced using the Lie-theoretic approach.

\begin{Example}
Let~\(n \ge 1\) and~\(m, k\) be nonzero integers.
\begin{enumerate}
\item The subgroup
\[\langle A_{1, n+1}^k A_{1,n}^m, A_{2, n+1}^k A_{2,n}^m, \ldots, A_{n-1, n+1}^k A_{n-1, n}^m, A_{n,n+1}^k \rangle\]
of~\(P_{n+1}\) is free if and only if~\((m, k) \ne (-1,-1)\).
\item The subgroup
\[\langle A_{1,n}^m A_{1, n+1}^k, A_{2,n}^m A_{2, n+1}^k, \ldots, A_{n-1, n}^m A_{n-1, n+1}^k, A_{n,n+1}^k \rangle\]
of~\(P_{n+1}\) is free if and only if~\((m, k) \ne (1, 1)\).
\end{enumerate}
\end{Example}

\subsection{Organization}
In Section~\ref{sec:main_results}, we construct a normal form for
a multiple HNN-extension of a free group by basis-conjugating embeddings (Theorem~\ref{t:diamond_lemma_for_G})
and give sufficient conditions to fulfill the requirements of the ping-pong lemma (Theorem~\ref{t:ping_pong_condition}).
In Section~\ref{sec:braid_group}, we recall what is known about freeness of subgroups in the pure braid group
and use results from Section~\ref{sec:main_results} to produce new families of free subgroups in the pure braid group (Theorem~\ref{t:main_theorem_braids}).
We have relegated to Appendix~\ref{sec:confluence_proof} the routine verification
of the local confluence of the rewriting system that appears in the proof of Theorem~\ref{t:diamond_lemma_for_G}.

\subsection*{Acknowledgment}
I am grateful to Leonid Danilevich for the statement and the proof of Lemma~\ref{t:danilevich} and to Artem Semidetnov for kindly introducing me to Stallings foldings.
I sincerely thank Andrei Malyutin for his careful reading of a draft of this paper.

\section{Proof of main results}
\label{sec:main_results}
Let~\(\Sigma\) be a set, and let~\(\Sigma^*\) denote the free monoid generated by~\(\Sigma\).
We call a relation on~\(\Sigma^*\) a \textit{congruence} if it is an equivalence relation such that~\(w_1uw_2 \sim w_1u'w_2\) whenever~\(u \sim u'\).
For~\(R \subset \Sigma^* \times \Sigma^*\), denote by~\(\rho(R)\) the congruence on~\(\Sigma^*\) generated by~\(R\). Let
\[M \defeq \mathrm{Mon}\langle \Sigma \mid R \rangle \overset{\mathrm{def}}{=} \Sigma^* / \rho(R)\]
be the monoid given by the presentation~\(\langle \Sigma \mid R \rangle\).
We have a natural quotient homomorphism
\[\pi \colon \Sigma^* \to M.\]
Since~\(R\) consists of ordered pairs, we can interpret it as a rewriting system:
each~\((r, r') \in R\) corresponds to the rewriting rule~\(w_1rw_2 \mapsto w_1r'w_2\).
We say that a rewriting system is \textit{locally confluent} if any two one-step applications of rules to the same word can be completed to chains of rule applications ending at a common word:
\begin{figure}[H]
\centering
\[\begin{tikzcd}[cramped]
	& w \\
	{w_1} && {w_2.} \\
	& {w'}
	\arrow["\forall"', from=1-2, to=2-1]
	\arrow["\forall", from=1-2, to=2-3]
	\arrow["\exists"', squiggly, from=2-1, to=3-2]
	\arrow["\exists", squiggly, from=2-3, to=3-2]
\end{tikzcd}\]
\caption{Here, a straight arrow corresponds to a single application of a rule, whereas a squiggly arrow corresponds to a finite sequence of applications.}
\end{figure}
We say that a rewriting system is \textit{Noetherian} if one cannot apply given rules to a word infinitely many times
(typically, this means that there is a bounded discrete parameter which is optimized after every application).
We say that a subset~\(N \subset \Sigma^*\) is a \textit{normal form} for~\(M\) if the restriction~\(\pi\mid_N \colon N \to M\) is a bijection.
A version of the diamond lemma (see \cite[Theorem~2]{N42}) says that if a rewriting system~\(\{r_j \mapsto r_j'\}\) is locally confluent and Noetherian, then the set~\(\{w \in \Sigma^* \mid w \not\supset r_j \; \forall j\}\) is a normal form for~\(M\).

The group~\(G\) given by the presentation~\eqref{eq:presentation_of_group} is generated as a monoid by the extended alphabet \(\Sigma = Y^\pm \cup X^\pm\).
Consider the following rewriting system for~\(G\):
\begin{enumerate}
\item~\(y_j^{\varepsilon} y_j^{-\varepsilon} \mapsto \varnothing\), where~\(1 \le j \le m\), and~\(\varepsilon \in \{\pm 1\}\);
\item~\(x_i^{\varepsilon} x_i^{-\varepsilon} \mapsto \varnothing\), where~\(1 \le i \le n\), and~\(\varepsilon \in \{\pm 1\}\);
\item~\(x_i v_{ij}^{-1} y_{ij}^\varepsilon \mapsto w_{ij}^{-1} y_{ij}^\varepsilon w_{ij} x_i v_{ij}^{-1}\), where~\(1 \le i \le n,\,1 \le j \le m_i\), and~\(\varepsilon \in \{\pm 1\}\);
\item~\(x_i^{-1} w_{ij}^{-1} y_{ij}^\varepsilon \mapsto v_{ij}^{-1} y_{ij}^\varepsilon v_{ij} x_i^{-1} w_{ij}^{-1}\), where~\(1 \le i \le n,\,1 \le j \le m_i\), and~\(\varepsilon \in \{\pm 1\}\).
\end{enumerate}

\begin{proof}[Proof of Theorem~\ref{t:diamond_lemma_for_G}]
Clearly, the above set (1)--(4) is a complete set of relations for~\(G\).
One can verify that this system is locally confluent.
The verification of this is technical and has been deferred to Appendix~\ref{sec:confluence_proof}.

In order to show that it is Noetherian, we introduce the following parameter.
Write a word~\(w \in \Sigma^*\) as a product
\[w = u_0 x_{i_1}^{\pm 1} u_1 x_{i_2}^{\pm 1} \ldots x_{i_k}^{\pm 1} u_k, \quad k \ge 0, \, u_i \in (Y^\pm)^*.\]
We can associate with the word~\(w\) the vector of non-negative numbers
\[\nu(w) = (|u_0|, \ldots, |u_k|) \in \ZZ_{\ge 0}^{k+1}.\]
Observe that each application of the rule either decreases~\(k\), or decreases one of the coordinates of~\(\nu(w)\) without touching the ones following it
(however, it can increase the previous ones).
Define an order on the set of tuples of non-negative integers by the rule~\((a_0, \ldots, a_k) < (b_0, \ldots, b_m)\) if~\(k < m\) or~\(k = m\) and for some~\(i\) one has~\(a_i < b_i\) and~\(a_j = b_j\) for~\(j > i\).
This order is known as the \textit{right shortlex order}.
The result follows, since this order is well-founded, hence there are no infinite decreasing chains.
\end{proof}

\begin{Remark}
\label{t:main_property_of_normal_form}
Let
\[w = u_0 x_{i_1}^{\varepsilon_1} u_1 x_{i_2}^{\varepsilon_2} \ldots x_{i_k}^{\varepsilon_k} u_k\]
be a word in~\(\Sigma^*\) and
\[w' = u_0' x_{j_1}^{\eta_1} u_1' x_{j_2}^{\eta_2} \ldots x_{j_l}^{\eta_l} u_l'\]
be its normal form,
where~\(u_i, u_j' \in (Y^\pm)^*\).
Then the sequence~\((x_{j_1}^{\eta_1}, \ldots, x_{j_l}^{\eta_l})\) is a subsequence of~\((x_{i_1}^{\varepsilon_1}, \ldots, x_{i_k}^{\varepsilon_k})\).
\end{Remark}

The following result is folklore.
\begin{Lemma}[{Ping-pong lemma, \cite[Lemma~4]{OL00}}]
Let~\(G\) be a group acting on a set~\(X\), and let~\(A_1, A_2, \ldots, A_k \subset G\) be nontrivial subgroups such that~\(|A_1| \ge 3\).
If there are disjoint nonempty subsets~\(U_1, \ldots, U_k \subset X\) such that for each pair of distinct indices~\(t \ne s\) and each nontrivial~\(a \in A_t\) one has~\(a(U_s) \subset U_t\), then~\(\langle A_1, \ldots, A_k \rangle = A_1 * \ldots * A_k\).
\end{Lemma}
Let~\(k \ge 1\), and let~\(T_1, \ldots, T_k\) be pairwise disjoint nonempty subsets of~\(X\).
Let~\(G\) act on itself by left shifts. For every~\(1 \le t \le k\) set
\begin{align*}
U_t = \{u \in G \mid & \; \text{the normal form of}\;u\;\text{looks like}\; u'x_i^{\pm 1} u'', \\
& \text{where}\;x_i\in T_t\;\text{and}\;u'\in (Y^\pm)^*\}.
\end{align*}
The sets~\(U_t\) are nonempty and pairwise disjoint.
We say that a subset~\(S\) of~\(G\) is \textit{greedy} if it possesses the following property:
if the concatenation of two words~\(uu'\) is already in the normal form, and~\(u \in S\), then~\(uu' \in S\).
Clearly, each~\(U_t\) is greedy.

\begin{proof}[Proof of Theorem~\ref{t:ping_pong_condition}]
Fix~\(1 \le t \ne s \le k\) and an element~\(a \in A_t \setminus \{1\}\).
We aim to show that~\(au \in U_t\) for each~\(u \in U_s\).

Let~\(u'x_i^\varepsilon u''\) be the normal form of~\(u\), where~\(x_i \in T_s\) and~\(u' \in (Y^\pm)^*\).
Let~\(p\) be the normal form of~\(au'\).
Since~\(a \in \langle Y, T_t \rangle\) and~\(u' \in \langle Y \rangle\), we have~\(p \in (Y^\pm \cup T_t^\pm)^*\) by Remark~\ref{t:main_property_of_normal_form}.
Hence the last letter of~\(p\) is not~\(x_i^{-\varepsilon}\), and the word~\(px_i^\varepsilon u''\) is the normal form of~\(au\).

The word~\(p\) contains~\(x_l^{\pm 1}\) for some~\(x_l \in T_t\), since if this is not the case, then~\(a = p {u'}^{-1} \in \langle Y \rangle\), which contradicts the assumptions.
Hence~\(p \in U_t\).
The result follows since~\(U_t\) is greedy. 
\end{proof}

Consider the special case where~\(\varphi_i = \psi_i\) for each~\(1 \le i \le n\).
In this case, we formulate a condition that ensures the second condition from Theorem~\ref{t:ping_pong_condition}.
The presentation~\(\eqref{eq:presentation_of_group}\) reads now as
\[G = \langle y_1, \ldots, y_m, x_1, \ldots, x_n \mid [x_i, y_{ij}^{w_{ij}}] = 1, \; 1 \le i \le n, \, 1 \le j \le m_i \rangle.\]
If we add the relations~\([x_i, y_j] = 1\) for all~\(i\) and~\(j\), then we obtain a map~\(\pi \colon G \surj F(X) \times F(Y)\).
Consider the following commutative diagram of projections:
\[\begin{tikzcd}
{F(X)} \\
{F(X \cup Y)} & G & {F(X) \times F(Y).} \\
{F(Y)}
\arrow["{\pi_X}", two heads, from=2-1, to=1-1]
\arrow["{[-]}", two heads, from=2-1, to=2-2]
\arrow["{\pi_Y}"', two heads, from=2-1, to=3-1]
\arrow["\pi", two heads, from=2-2, to=2-3]
\arrow["{\pi_X}"', two heads, from=2-3, to=1-1]
\arrow["{\pi_Y}", two heads, from=2-3, to=3-1]
\end{tikzcd}\]
\begin{Proposition}
\label{t:orbit_set_is_free_basis}
Let~\(\varphi\) be an automorphism of a free group~\(F(X \cup Y)\) that descends as the identity map on~\(F(X) \times F(Y)\), let~\(w \in F(X \cup Y)\) be an element with~\(\pi_X(w) \ne 1\),
and let~\(A\) be a subgroup of~\(G\) generated by the set~\(\{[\varphi^k(w)] \mid k \in \ZZ\}\). Then the intersection~\(A \cap \langle Y \rangle\) is trivial.
\end{Proposition}
\begin{proof}
Note that the restriction~\([-]\mid_{F(Y)} \colon F(Y) \to G\) is injective, since~\(G\) is an HNN-extension of~\(F(Y)\).
Hence the restriction
\[\pi \mid_{\langle Y \rangle} \colon \langle Y \rangle \to F(X) \times F(Y)\]
is injective. Therefore, it is enough to show that
\[\pi(A) \cap \pi(\langle Y \rangle) = 1.\]
Note that
\begin{align*}
\pi(A) &= \pi(\langle [\varphi^k(w)] \mid k \in \ZZ \rangle) = \langle (\pi_X(\varphi^k(w)), \pi_Y(\varphi^k(w))) \mid k \in \ZZ \rangle \\
&= \langle (\pi_X(w), \pi_Y(w)) \rangle
\end{align*}
and~\(\pi(\langle Y \rangle) = 1 \times F(Y)\).
The result follows since the restriction of~\(\pi_X\) to~\(\langle (\pi_X(w), \pi_Y(w))\rangle\) is injective
while~\(\pi_X(1 \times F(Y)) = 1\).
\end{proof}

\section{Free subgroups in pure braid group}
\label{sec:braid_group}
We now apply the general machinery developed above to the pure braid group.
The pure braid group contains many free subgroups.
We list some of them below.
\begin{enumerate}
\item For every~\(1 \le i \le n+1\) there is a strand-removing homomorphism~\(d_i \colon P_{n+1} \to P_n\). It is well known that its kernel~\(\ker(d_i)\) is free of rank~\(n\). The proof is homotopy-theoretic.
\item There is a rank~\(n\) free subgroup of~\(P_{n+1}\) given by an embedding of Milnor's construction~\(F[S^1]\) into a simplicial group~\(\AP\) of pure braids. The classical proof \cite{CW04, CW11} relies on the Lie rings~\(\mathrm{gr}_*(F_n)\) and~\(\mathrm{gr}_*(P_{n+1})\). There is also a group-theoretic proof \cite{Io25}.
\item The images of the above subgroups under automorphisms of the pure braid group are also free.
These include, for example, the free subgroup of rank~\(n\) generated by the braids\footnote{See \cite[Fig.~1]{Io25} for a picture of the braid~\(A_{0,i}\).}~\(A_{0,1}, \ldots, A_{0,n}\), where
\[A_{0,i} = (A_{1,i} A_{2,i} \ldots A_{i-1, i} A_{i, i+1} \ldots A_{i,n+1})^{-1}.\]
Note that~\(\langle A_{0, 1}, \ldots, A_{0, n+1} \rangle\) is a subgroup of braids that become trivial after being transplanted
onto a two-sphere.
It is worth noting that by the main result of \cite{DL16} we have an isomorphism
\[\ker(P_{n+1} \to P_{n+1}(S^2)) \cong F_n \times \ZZ.\]
\item In \cite{LM10}, it is proved that each pair of non-commuting braids~\(x\) and~\(y\) in \(P_{n+1}\) generate a free subgroup of rank~\(2\).
The proof involves the theory of~\(3\)-manifolds and the theory of group actions on trees.
This result suggests that there might be a lot of free subgroups in general.
\end{enumerate}

We refer to \cite{MK99} for standard notions of braid theory.
The group~\(P_{n+1}^{(2)}\) decomposes as a semidirect product
\begin{equation}\label{splitting_for_p2}
P_{n+1}^{(2)} = F(x_1, \ldots, x_n) \rtimes F(y_1, \ldots, y_{n-1}),
\end{equation}
and the action of~\(y_j = A_{j,n}\) on~\(x_i = A_{i,n+1}\) is given by the following relations:
\begin{tagcases}[x_i^{y_j} =]
&x_i, && \quad i < j; \label{rel1}\tag{R1} \\
&x_i^{(x_i x_n)^{-1}}, && \quad i = j; \label{rel2}\tag{R2} \\
&x_i^{[x_n, x_j]}, && \quad j < i < n; \label{rel3}\tag{R3} \\
&x_n^{x_j^{-1}}, && \quad j < i = n. \label{rel4}\tag{R4}
\end{tagcases}

We start with some manipulations with relations in the group~\(P_{n+1}^{(2)}\).
\begin{itemize}
\item The relation~\eqref{rel1} reads as
\begin{equation}
[x_i, y_j] = 1. \label{rel1new}\tag{R1'}
\end{equation}
\item The relation~\eqref{rel2} reads as~\(x_i^{y_i x_i x_n} = x_i\), which is equivalent by~\eqref{rel4} to
\begin{equation}
x_i^{x_n} = x_i^{y_i^{-1}}. \label{rel2new}\tag{R2'}
\end{equation}
\item The relation~\eqref{rel3} reads as~\(x_i^{y_j x_j x_n} = x_i^{x_n x_j}\), which is equivalent to~\(x_i^{x_n} = x_i^{x_n y_j^{-1}}\) by~\eqref{rel4}, which is equivalent to~\(x_i^{x_n} = x_i^{y_i^{-1} y_j}\) by~\eqref{rel2new}. By~\eqref{rel4}, one can replace this relation with
\begin{equation}
[x_i, y_j^{y_i}] = 1. \label{rel3new}\tag{R3'}
\end{equation}
\item By~\eqref{rel2new}, the relation~\eqref{rel4} is equivalent to
\begin{equation}
y_i^{x_n} = y_i [x_i, y_i]. \label{rel4new}\tag{R4'}
\end{equation} 
\end{itemize}
Hence
\begin{align*}
P_{n+1}^{(2)} &= \langle x_i, y_i \mid \text{(\ref{rel1})}, \text{(\ref{rel2})}, \text{(\ref{rel3})}, \text{(\ref{rel4})} \rangle \\
&= \langle x_i, y_i \mid \text{(\ref{rel1new})}, \text{(\ref{rel2new})}, \text{(\ref{rel3new})}, \text{(\ref{rel4new})} \rangle.
\end{align*}

Let~\(G_n\) be the group given by the following presentation
\[G_n \defeq \langle x_1, \ldots, x_{n-1}, y_1, \ldots, y_{n-1} \mid [x_i, y_j] = 1, \, i < j; \; [x_i, y_j^{y_i}] = 1, \, i > j \rangle.\]
The group~\(P_{n+1}^{(2)}\) admits a semidirect product decomposition
\[P_{n+1}^{(2)} = G_n \rtimes \langle t \rangle,\]
where the letter~\(t \defeq x_n\) acts on~\(G_n\) via the automorphism~\(\varphi \colon G_n \to G_n\) given by
\begin{equation}
\label{eq:phi_formula}
\varphi(x_i) = x_i^{y_i^{-1}}, \; \varphi(y_i) = y_i [x_i, y_i].
\end{equation}

Let~\(X \defeq \{x_1, \ldots, x_{n-1}\}\) and~\(Y \defeq \{y_1, \ldots, y_{n-1}\}\).
\begin{Corollary}[Normal form for~\(G_n\)]
Each element of~\(G_n\)
can be uniquely written as a reduced word~\(w\) in the extended alphabet~\(X^\pm \cup Y^\pm\) such that~\(w\)
does not contain subwords~\(x_i^\varepsilon y_j^\eta\) and~\(x_j^\varepsilon y_j^{-1} y_i^\eta\) with~\(i < j\) and~\(\varepsilon, \eta \in \{\pm 1\}\).
\end{Corollary}

\begin{Example}
The group~\(P_4\) decomposes as the semidirect product
\[F(x_1, x_2, x_3) \rtimes F(y_1, y_2) \rtimes \langle A_{1,2} \rangle.\]
If we replace~\(A_{1,2}\) with the generator of the center
\[z = (A_{1,2} A_{1,3} A_{1,4}) \cdot (A_{2,3} A_{2,4}) \cdot A_{3,4} = (A_{1,2} y_1 x_1) \cdot (y_2 x_2) \cdot x_3,\]
then we can rewrite this decomposition as
\begin{align*}
P_4 &= P^{(2)}_4 \oplus \ZZ = (G_3 \rtimes \ZZ) \oplus \ZZ \\
&= \langle y_1, y_2, x_1, x_2 \mid [x_1, y_2] = 1, [x_2, y_1^{y_2}] = 1 \rangle \rtimes (\ZZ \oplus \ZZ),
\end{align*}
where the generators of~\(\ZZ \oplus \ZZ\) are~\(x_3\) and~\(z\).
The letter~\(x_3\) acts on the normal subgroup~\(\langle y_1, y_2, x_1, x_2 \rangle\) of~\(P_4\) by the automorphism~\(\varphi\) given by the formula~\eqref{eq:phi_formula}, and~\(z\) is central.
\end{Example}

\begin{Lemma}
\label{t:danilevich}
Let~\(G\) be a group,~\(\varphi\) an automorphism of~\(G\), and~\(H\) a subgroup of~\(G\).
Consider the semidirect product~\(G \rtimes_\varphi \ZZ\) of~\(G\) with the infinite cyclic group~\(\ZZ = \langle t \rangle\).
The following are equivalent\textup{:}
\begin{enumerate}
\item the natural map~\(\coprod_{k \in \ZZ} \varphi^{k}(H)\to G\) is injective\textup{;}
\item the subgroup~\(\langle H, t \rangle\) of~\(G \rtimes_\varphi \ZZ\) decomposes as the free product~\(H \ast \langle t \rangle\).
\end{enumerate}
\end{Lemma}
\begin{proof}
\((1) \Rightarrow (2)\). Let~\(\Phi \colon H * \ZZ \to G \rtimes \ZZ\) be the natural map.
Let
\[w = t^{n_1} h_1 t^{n_2} h_2 \ldots h_{k-1} t^{n_k}\]
be an element in the free product~\(H * \ZZ\) written in the normal form, i.e., we assume that~\(n_i \ne 0\) for all~\(i \notin \{1, k\}\), and~\(h_i \ne 1\) for all~\(i\).
One has
\begin{align*}
\Phi(w) &= t^{n_1} h_1 t^{-n_1} t^{n_1+n_2} h_2 t^{-(n_1+n_2)} \ldots h_{k-1} t^{-(n_1+\ldots+n_{k-1})} \cdot t^{n_1+n_2+\ldots+n_k} \\
&= \varphi^{-n_1}(h_1) \varphi^{-(n_1+n_2)}(h_2) \ldots \varphi^{-(n_1+\ldots+n_{k-1})}(h_{k-1}) \cdot t^{n_1+\ldots+n_k}.
\end{align*}
In the above expansion, all adjacent powers of~\(\varphi\) are distinct since~\(n_i \ne 0\).
Hence if~\(k > 1\), then~\(\Phi(w) \ne 1\). If~\(k = 1\), then~\(\Phi(w) = t^{n_1}\), which is trivial if and only if~\(n_1 = 0\), i.e., if~\(w = 1\).

\((2) \Rightarrow (1)\). This follows since we have the natural embedding
\[\coprod_{k \in \ZZ} t^{-k} H t^k \hookrightarrow H * \langle t \rangle.\qedhere\]
\end{proof}

\begin{proof}[Proof of Theorem~\ref{t:main_theorem_braids}]
Let~\(B_i = \langle w_i \rangle\) for~\(1 \le i < n\). We aim to show that the natural map
\[B_1 * \ldots * B_{n-1} * \langle t \rangle \to G_n \rtimes \langle t \rangle = P_{n+1}^{(2)}\]
is injective.

Let~\(A_i = \langle \varphi^k(w_i) \mid k \in \ZZ \rangle\) for~\(1 \le i < n\).
In the light of Lemma~\ref{t:danilevich}, it is enough to show that:
\begin{enumerate}
\item for each~\(1 \le i < n\), the set~\(\{\varphi^k(w_i)\}_{k \in \ZZ}\) is a free basis of some subgroup in~\(G_n\);
\item the natural map~\(A_1 * \ldots * A_{n-1} \to G_n\) is injective.
\end{enumerate}
By the main result of \cite{LM10}, the condition~\([w_i, t] \ne 1\) implies that the subgroup~\(\langle w_i, t \rangle\) is free of rank two.
The other direction of Lemma~\ref{t:danilevich} implies the statement (1).

The formula~\eqref{eq:phi_formula} lifts the automorphism~\(\varphi \colon G_n \to G_n\)
to an automorphism of~\(F(X \cup Y)\).
The latter induces the identity map on the product~\(F(X) \times F(Y)\). In~\(G_n\) we have~\(A_i \cap \langle Y \rangle = 1\) by Proposition~\ref{t:orbit_set_is_free_basis}. Then, the statement (2) follows from Theorem~\ref{t:ping_pong_condition}.
\end{proof}

\subsection{Possible generalization}
It is natural to ask whether the statement of Theorem~\ref{t:main_theorem_braids} still hold true if the last generator~\(t = x_n\) is replaced by a word~\(w_n\) such that
\[w_n\in \langle Y, x_n \rangle, \; \exp_{x_n} (w_n) \ne 0, \; \text{and} \; [w_i, w_n] \ne 1 \; \text{for every} \; i < n.\]
It seems that in this case, the subgroup~\(\langle w_1, \ldots, w_n \rangle\) would also be free, but the argument has to be modified: the presented proof uses the fact that the conjugation by the last generator induces the identity automorphism on~\(F(X) \times F(Y)\)
in order to demonstrate that~\(\langle w_n^{-k} w_i w_n^k \mid k \in \ZZ \rangle \cap \langle Y \rangle = 1\).
Our argument easily generalizes to words~\(w_n\) from the normal closure of~\(x_n\) in the group~\(\langle Y, x_n \rangle\).
However, for the case of three strands we can ensure that no restrictions are needed.
Note that~\(G_2 = \langle y_1, x_1 \rangle\) and~\(P^{(2)}_3 = P_3\).
\begin{Proposition}
Let~\(w_1, w_2 \in P_3\) be braids such that~\(w_1 \in \langle y_1, x_1 \rangle\) and~\(\exp_{x_1}(w_1) \ne 0\).
Then the intersection~\(\langle w_2^{-k} w_1 w_2^k \mid k \in \ZZ\rangle \cap \langle y_1 \rangle\)
is trivial.
\end{Proposition}
\begin{proof}
The center of~\(P_3\) is generated by the braid~\(z = y_1 x_1 t\), where~\(t = x_2\) as before. Hence we can assume that~\(w_2\) is from~\(\langle y_1, x_1 \rangle\).
It is sufficient to demonstrate that the normal closure of~\(w_1\) in~\(\langle y_1, x_1 \rangle\) does not intersect~\(\langle y_1 \rangle\).
If~\(w_1 = (u^k)^v\), then the normal closure of~\(w_1\) is contained in the normal closure of~\(u\). Hence we can assume that the cyclic reduction of~\(w_1\) is not a proper power.

Note that~\(\langle y_1, x_1 \rangle\) is free of rank two. The one-relator group
\[\langle x_1, y_1 \mid w_1 \rangle\]
is torsion-free by \cite[Theorem~4.12]{MKS04}. It follows that any nontrivial intersection should give~\(y_1 = 1\) in this group. Hence this group would be generated by~\(x_1\). But since~\(\exp_{x_1}(w_1) \ne 0\), this would imply that this group is finite cyclic, which is impossible.
\end{proof}

\subsection{Lie-theoretic approach}
With a group~\(G\) filtered by the lower central series one can associate a graded Lie ring
\[\mathrm{gr}_*(G) \defeq \bigoplus_{k \ge 1} \frac{\gamma_k(G)}{\gamma_{k+1}(G)}\]
with the Lie bracket given by the commutator.
This construction is functorial.
Moreover, if~\(G\) is residually nilpotent, then injectivity of a homomorphism~\(f \colon G \to H\)
follows from injectivity of the homomorphism~\(\mathrm{gr}_*(f)\).
This approach was taken in \cite{CW04, CW11} to show injectivity of the natural simplicial map
\(F[S^1] \to \AP\).
In degree~\(n\), it corresponds to the homomorphism~\(F_n \to P_{n+1}\).

Note that if two homomorphisms~\(f, g\colon G \to H\) induce the same homomorphism of abelianizations~\(G_\mathrm{ab} \to H_{\mathrm{ab}}\),
then they induce the same homomorphism of graded Lie rings~\(\mathrm{gr}_*(G) \to \mathrm{gr}_*(H)\).
Since the subgroup
\[\langle y_1 x_1, \ldots, y_{n-1} x_{n-1}, t \rangle \subset P_{n+1}\]
is isomorphic to~\(F_{n-1} \times \ZZ\), it follows that the corresponding homomorphism~\(F_n \to P_{n+1}\) is not injective.
Therefore, the induced homomorphism~\(\mathrm{gr}_*(F_n) \to \mathrm{gr}_*(P_{n+1})\) is not injective.
On the other hand, Theorem~\ref{t:main_theorem_braids} implies that the subgroup
\[\langle x_1 y_1, \ldots, x_{n-1} y_{n-1}, t\rangle \subset P_{n+1}\]
is free of rank~\(n\).
But this clearly could not be deduced via Lie-theoretic approach, since~\(y_ix_i = [y_i, x_i] x_iy_i\) and the two homomorphisms~\(F_n \to P_{n+1}\) induce the same map of abelianizations.

\appendix
\section{Proof of local confluence}
\label{sec:confluence_proof}
In this section, we fill a gap in the proof of Theorem~\ref{t:diamond_lemma_for_G} and show that
the rewriting system is locally confluent.

In order to demonstrate that a rewriting system is locally confluent, one must inspect all possible ways to apply
two rules in some overlapping fragments of a given word.
It is easy to see that the left hand sides of the rules~(3) and~(4) do not intersect
and there are no self-intersections for the rules~(3) and~(4).
The case of self-overlappings for the rules~(1) and~(2) is analysed in Fig.~\ref{fig:1_2}.
The cases of overlappings of the rules~(1) and~(3), and of the rules~(2) and~(3), are analysed
in Fig.~\ref{fig:1_3} and Fig.~\ref{fig:2_3} respectively.
The cases of overlappings of the rules~(1) and~(4), and of the rules~(2) and~(4), are completely analogous.

\tikzset{every picture/.style={line width=0.75pt, scale=0.72}}

\begin{figure}[H]
\centering
\begin{tikzpicture}[x=0.75pt,y=0.75pt,yscale=-1,xscale=1]
\draw [color=darkgreen, draw opacity=1] (380,40) -- (398.89,68.34);
\draw [shift={(400,70)}, rotate = 236.31] [color=darkgreen, draw opacity=1][line width=0.75] (10.93,-3.29) .. controls (6.95,-1.4) and (3.31,-0.3) .. (0,0) .. controls (3.31,0.3) and (6.95,1.4) .. (10.93,3.29);
\draw [color=blue, draw opacity=1] (320,40) -- (301.11,68.34);
\draw [shift={(300,70)}, rotate = 303.69] [color=blue, draw opacity=1][line width=0.75] (10.93,-3.29) .. controls (6.95,-1.4) and (3.31,-0.3) .. (0,0) .. controls (3.31,0.3) and (6.95,1.4) .. (10.93,3.29);
\draw [color=darkgreen, draw opacity=1] (170,40) -- (188.89,68.34);
\draw [shift={(190,70)}, rotate = 236.31] [color=darkgreen, draw opacity=1][line width=0.75] (10.93,-3.29) .. controls (6.95,-1.4) and (3.31,-0.3) .. (0,0) .. controls (3.31,0.3) and (6.95,1.4) .. (10.93,3.29);
\draw [color=blue, draw opacity=1] (110,40) -- (91.11,68.34);
\draw [shift={(90,70)}, rotate = 303.69] [color=blue, draw opacity=1][line width=0.75] (10.93,-3.29) .. controls (6.95,-1.4) and (3.31,-0.3) .. (0,0) .. controls (3.31,0.3) and (6.95,1.4) .. (10.93,3.29);

\draw (89,9) node [anchor=north west][inner sep=0.75pt] {$\dotsc \uc{$x_{i}^{\varepsilon}$}{blue} \uc[1]{\uc{$x_{i}^{-\varepsilon}$}{blue}}{darkgreen} \uc[2]{$x_{i}^{\varepsilon}$}{darkgreen} \dotsc $};
\draw (49,69) node [anchor=north west][inner sep=0.75pt]  {$\dotsc x_{i}^{\varepsilon} \dotsc $};
\draw (169,69) node [anchor=north west][inner sep=0.75pt]  {$\dotsc x_{i}^{\varepsilon} \dotsc $};
\draw (299,9) node [anchor=north west][inner sep=0.75pt]  {$\dotsc \uc[1.12]{$y_{j}^{\varepsilon}$}{blue} \uc[1]{\uc{$y_{j}^{-\varepsilon}$}{blue}}{darkgreen} \uc[2.12]{$y_{j}^{\varepsilon}$}{darkgreen} \dotsc $};
\draw (259,69) node [anchor=north west][inner sep=0.75pt]  {$\dotsc y_{j}^{\varepsilon} \dotsc $};
\draw (389,69) node [anchor=north west][inner sep=0.75pt]  {$\dotsc y_{j}^{\varepsilon} \dotsc $};
\draw (398,39) node [anchor=north west][inner sep=0.75pt] [align=left] {\textcolor{darkgreen}{(1)}};
\draw (189,39) node [anchor=north west][inner sep=0.75pt] [align=left] {\textcolor{darkgreen}{(2)}};
\draw (69,39) node [anchor=north west][inner sep=0.75pt] [align=left] {\textcolor{blue}{(2)}};
\draw (279,39) node [anchor=north west][inner sep=0.75pt] [align=left] {\textcolor{blue}{(1)}};
\end{tikzpicture}
\caption{Self-overlapping for rules (1) and (2).}
\label{fig:1_2}
\end{figure}

\begin{figure}[H]
\centering
\begin{tikzpicture}[x=0.75pt,y=0.75pt,yscale=-1,xscale=1]
\draw [color=blue, draw opacity=1] (110,530) -- (91.11,558.34);
\draw [shift={(90,560)}, rotate = 303.69] [color=blue, draw opacity=1][line width=0.75] (10.93,-3.29) .. controls (6.95,-1.4) and (3.31,-0.3) .. (0,0) .. controls (3.31,0.3) and (6.95,1.4) .. (10.93,3.29);
\draw [color=darkgreen, draw opacity=1] (180,530) -- (198.89,558.34);
\draw [shift={(200,560)}, rotate = 236.31] [color=darkgreen, draw opacity=1][line width=0.75] (10.93,-3.29) .. controls (6.95,-1.4) and (3.31,-0.3) .. (0,0) .. controls (3.31,0.3) and (6.95,1.4) .. (10.93,3.29);
\draw [color=darkgreen, draw opacity=1] (260,590) -- (260,618);
\draw [shift={(260,620)}, rotate = 270] [color=darkgreen, draw opacity=1][line width=0.75] (10.93,-3.29) .. controls (6.95,-1.4) and (3.31,-0.3) .. (0,0) .. controls (3.31,0.3) and (6.95,1.4) .. (10.93,3.29);
\draw [color=darkgreen, draw opacity=1] (260,650) .. controls (261.67,651.67) and (261.67,653.33) .. (260,655) .. controls (258.33,656.67) and (258.33,658.33) .. (260,660) .. controls (261.67,661.67) and (261.67,663.33) .. (260,665) .. controls (258.33,666.67) and (258.33,668.33) .. (260,670) -- (260,678);
\draw [shift={(260,680)}, rotate = 270] [color=darkgreen, draw opacity=1][line width=0.75] (10.93,-3.29) .. controls (6.95,-1.4) and (3.31,-0.3) .. (0,0) .. controls (3.31,0.3) and (6.95,1.4) .. (10.93,3.29);
\draw [color=darkgreen, draw opacity=1] (260,710) -- (260,738);
\draw [shift={(260,740)}, rotate = 270] [color=darkgreen, draw opacity=1][line width=0.75] (10.93,-3.29) .. controls (6.95,-1.4) and (3.31,-0.3) .. (0,0) .. controls (3.31,0.3) and (6.95,1.4) .. (10.93,3.29);
\draw [color=darkgreen, draw opacity=1] (260,770) .. controls (261.67,771.67) and (261.67,773.33) .. (260,775) .. controls (258.33,776.67) and (258.33,778.33) .. (260,780) .. controls (261.67,781.67) and (261.67,783.33) .. (260,785) .. controls (258.33,786.67) and (258.33,788.33) .. (260,790) -- (260,798);
\draw [shift={(260,800)}, rotate = 270] [color=darkgreen, draw opacity=1][line width=0.75] (10.93,-3.29) .. controls (6.95,-1.4) and (3.31,-0.3) .. (0,0) .. controls (3.31,0.3) and (6.95,1.4) .. (10.93,3.29);

\draw (79,499) node [anchor=north west][inner sep=0.75pt]  {$\dotsc \uc{$x_{i} v_{ij}^{-1}$}{darkgreen} \uc[1]{\uc[1.12]{$y_{ij}^{\varepsilon}$}{darkgreen}}{blue} \uc[2]{$y_{ij}^{-\varepsilon}$}{blue} \dotsc $};
\draw (39,559) node [anchor=north west][inner sep=0.75pt]  {$\dotsc x_{i} v_{ij}^{-1} \dotsc $};
\draw (169,559) node [anchor=north west][inner sep=0.75pt]  {$\dotsc w_{ij}^{-1} y_{ij}^{\varepsilon} w_{ij} \uc{$x_{i} v_{ij}^{-1} y_{ij}^{-\varepsilon}$}{darkgreen} \dotsc $};
\draw (149,619) node [anchor=north west][inner sep=0.75pt]  {$\dotsc w_{ij}^{-1} y_{ij}^{\varepsilon} \uc{$w_{ij} w_{ij}^{-1}$}{darkgreen} y_{ij}^{-\varepsilon} w_{ij} x_{i} v_{ij}^{-1} \dotsc $};
\draw (169,679) node [anchor=north west][inner sep=0.75pt]  {$\dotsc w_{ij}^{-1} \uc{$y_{ij}^{\varepsilon} y_{ij}^{-\varepsilon}$}{darkgreen} w_{ij} x_{i} v_{ij}^{-1} \dotsc $};
\draw (199,739) node [anchor=north west][inner sep=0.75pt]  {$\dotsc \uc{$w_{ij}^{-1} w_{ij}$}{darkgreen} x_{i} v_{ij}^{-1} \dotsc $};
\draw (221,799) node [anchor=north west][inner sep=0.75pt]  {$\dotsc x_{i} v_{ij}^{-1} \dotsc $};
\draw (69,529) node [anchor=north west][inner sep=0.75pt] [align=left] {\textcolor{blue}{(1)}};
\draw (199,529) node [anchor=north west][inner sep=0.75pt] [align=left] {\textcolor{darkgreen}{(3)}};
\draw (269,589) node [anchor=north west][inner sep=0.75pt] [align=left] {\textcolor{darkgreen}{(3)}};
\draw (269,649) node [anchor=north west][inner sep=0.75pt] [align=left] {\textcolor{darkgreen}{(1)}};
\draw (269,709) node [anchor=north west][inner sep=0.75pt] [align=left] {\textcolor{darkgreen}{(1)}};
\draw (269,769) node [anchor=north west][inner sep=0.75pt] [align=left] {\textcolor{darkgreen}{(1)}};
\end{tikzpicture}
\caption{Overlapping of rules (1) and (3).}
\label{fig:1_3}
\end{figure}

\begin{figure}[H]
\centering
\begin{tikzpicture}[x=0.75pt,y=0.75pt,yscale=-1,xscale=1]
\draw [color=blue, draw opacity=1] (110,170) -- (91.11,198.34);
\draw [shift={(90,200)}, rotate = 303.69] [color=blue, draw opacity=1][line width=0.75] (10.93,-3.29) .. controls (6.95,-1.4) and (3.31,-0.3) .. (0,0) .. controls (3.31,0.3) and (6.95,1.4) .. (10.93,3.29);
\draw [color=darkgreen, draw opacity=1] (180,170) -- (198.89,198.34);
\draw [shift={(200,200)}, rotate = 236.31] [color=darkgreen, draw opacity=1][line width=0.75] (10.93,-3.29) .. controls (6.95,-1.4) and (3.31,-0.3) .. (0,0) .. controls (3.31,0.3) and (6.95,1.4) .. (10.93,3.29);
\draw [color=darkgreen, draw opacity=1] (260,230) -- (260,258);
\draw [shift={(260,260)}, rotate = 270] [color=darkgreen, draw opacity=1][line width=0.75] (10.93,-3.29) .. controls (6.95,-1.4) and (3.31,-0.3) .. (0,0) .. controls (3.31,0.3) and (6.95,1.4) .. (10.93,3.29);
\draw [color=darkgreen, draw opacity=1] (260,290) .. controls (261.67,291.67) and (261.67,293.33) .. (260,295) .. controls (258.33,296.67) and (258.33,298.33) .. (260,300) .. controls (261.67,301.67) and (261.67,303.33) .. (260,305) .. controls (258.33,306.67) and (258.33,308.33) .. (260,310) -- (260,318);
\draw [shift={(260,320)}, rotate = 270] [color=darkgreen, draw opacity=1][line width=0.75] (10.93,-3.29) .. controls (6.95,-1.4) and (3.31,-0.3) .. (0,0) .. controls (3.31,0.3) and (6.95,1.4) .. (10.93,3.29);
\draw [color=darkgreen, draw opacity=1] (260,350) -- (260,378);
\draw [shift={(260,380)}, rotate = 270] [color=darkgreen, draw opacity=1][line width=0.75] (10.93,-3.29) .. controls (6.95,-1.4) and (3.31,-0.3) .. (0,0) .. controls (3.31,0.3) and (6.95,1.4) .. (10.93,3.29);
\draw [color=darkgreen, draw opacity=1] (260,410) .. controls (261.67,411.67) and (261.67,413.33) .. (260,415) .. controls (258.33,416.67) and (258.33,418.33) .. (260,420) .. controls (261.67,421.67) and (261.67,423.33) .. (260,425) .. controls (258.33,426.67) and (258.33,428.33) .. (260,430) -- (260,438);
\draw [shift={(260,440)}, rotate = 270] [color=darkgreen, draw opacity=1][line width=0.75] (10.93,-3.29) .. controls (6.95,-1.4) and (3.31,-0.3) .. (0,0) .. controls (3.31,0.3) and (6.95,1.4) .. (10.93,3.29);

\draw (80,139) node [anchor=north west][inner sep=0.75pt]  {$\dotsc \uc{$x_{i}^{-1}$}{blue} \uc[1]{\uc{$x_{i}$}{blue}}{darkgreen} \uc[1.68]{$v_{ij}^{-1} y_{ij}^{\varepsilon}$}{darkgreen} \dotsc $};
\draw (169,199) node [anchor=north west][inner sep=0.75pt]  {$\dotsc \uc{$x_{i}^{-1} w_{ij}^{-1} y_{ij}^{\varepsilon}$}{darkgreen} w_{ij} x_{i} v_{ij}^{-1} \dotsc $};
\draw (39,199) node [anchor=north west][inner sep=0.75pt]  {$\dotsc v_{ij}^{-1} y_{ij}^{\varepsilon} \dotsc $};
\draw (159,259) node [anchor=north west][inner sep=0.75pt]  {$\dotsc v_{ij}^{-1} y_{ij}^{\varepsilon} v_{ij} x_{i}^{-1} \uc{$w_{ij}^{-1} w_{ij}$}{darkgreen} x_{i} v_{ij}^{-1} \dotsc $};
\draw (179,319) node [anchor=north west][inner sep=0.75pt]  {$\dotsc v_{ij}^{-1} y_{ij}^{\varepsilon} v_{ij} \uc{$x_{i}^{-1} x_{i}$}{darkgreen} v_{ij}^{-1} \dotsc $};
\draw (199,379) node [anchor=north west][inner sep=0.75pt]  {$\dotsc v_{ij}^{-1} y_{ij}^{\varepsilon} \uc{$v_{ij} v_{ij}^{-1}$}{darkgreen} \dotsc $};
\draw (219,439) node [anchor=north west][inner sep=0.75pt]  {$\dotsc v_{ij}^{-1} y_{ij}^{\varepsilon} \dotsc $};
\draw (69,169) node [anchor=north west][inner sep=0.75pt] [align=left] {\textcolor{blue}{(2)}};
\draw (199,169) node [anchor=north west][inner sep=0.75pt] [align=left] {\textcolor{darkgreen}{(3)}};
\draw (269,229) node [anchor=north west][inner sep=0.75pt] [align=left] {\textcolor{darkgreen}{(4)}};
\draw (269,289) node [anchor=north west][inner sep=0.75pt] [align=left] {\textcolor{darkgreen}{(1)}};
\draw (269,349) node [anchor=north west][inner sep=0.75pt] [align=left] {\textcolor{darkgreen}{(2)}};
\draw (269,409) node [anchor=north west][inner sep=0.75pt] [align=left] {\textcolor{darkgreen}{(1)}};
\end{tikzpicture}
\caption{Overlapping of rules (2) and (3).}
\label{fig:2_3}
\end{figure}

\end{document}